\documentclass[final,11pt]{article}
\usepackage[cp1251]{inputenc}

\newtheorem{theorem}{Theorem}[section]
\newtheorem{lemma}[theorem]{Lemma}

\newtheorem{definition}[theorem]{Definition}

\usepackage{amsmath}
    \usepackage{bm}
\usepackage{amsfonts,amssymb}

  \usepackage{xcolor}
\usepackage{epsfig}
\usepackage{psfrag}
 \usepackage{epsfrag}

\newtheorem{remark}[theorem]{Remark}

\newcommand{\bfalf}{\mbox{\boldmath $\alpha$}}

\def\eqsp{\noalign{\vskip 5pt}}
\def\dsp{\displaystyle}
\newcommand{\bftau}{\mbox{\boldmath $\tau$}}
\newcommand{\bfxi}{\mbox{\boldmath $\xi$}}

\def\dsp{\displaystyle}
\newcommand{\eqml}[1]{\eql{#1}\begin{array}{rcl}}
\newcommand{\enml}{\end{array}\end{equation}}
\newcommand{\eqm}{\begin{eqnarray}}
\newcommand{\enm}{\end{eqnarray}}
\newcommand{\eqmno}{\begin{eqnarray*}}
\newcommand{\enmno}{\end{eqnarray*}}
\newcommand{\goto}{\rightarrow}
\newcommand{\eql}[1]{\begin{equation}\label{#1}}
\newcommand{\half}{{\frac{1}{2}}}
\newcommand{\lp}{\left (}
\newcommand{\rp}{\right )}
\newcommand{\vthin}{\vskip 5pt}
\newcommand{\eqn}[1]{(\ref{#1})}
\def\grad{\nabla}
\def\qquad{\quad\quad}
\def\bi{\begin{itemize}}
\def\ei{\end{itemize}}
\def\be{\begin{enumerate}}
\def\ee{\end{enumerate}}
\newcommand{\ignore}[1]{}

\begin{document}
\title{Convergence study of   IB methods for Stokes equations with no-slip boundary conditions} 

\author{Zhilin Li \thanks{
Department of Mathematics, North Carolina State University, Raleigh,
NC 27695-8205, USA, Email: zhilin@math.ncsu.edu}
\and Kejia Pan\thanks{School of Mathematics and Statistics, HNP-LAMA, Central South University, Changsha  410083, Hunan, China,  Email: kejiapan@csu.edu.cn}
\and Juan Ruiz-\'Alvarez \thanks{Departamento de Matem\'atica  Aplicada y Estad\'istica. Universidad  Polit\'ecnica de Cartagena,  Spain, Email: juan.ruiz@upct.es}
}


\maketitle

\markright{Convergence of IB method}

\begin{abstract}
Peskin's Immersed Boundary (IB) model and method  are  among  one of the most important
modeling tools and numerical methods. The IB method has been known to be first
order accurate in the velocity.  However, almost no  rigorous theoretical proof can be found in the literature for
Stokes  equations with a prescribed velocity  boundary condition.
In this paper,  it has been shown that the pressure of the Stokes equation has a convergence order
$O(\sqrt{h} |\log h| )$  in the $L^2$ norm while the velocity has an
$O(h |\log h| )$ convergence order in the infinity norm in two-space dimensions. The proofs are based on 
splitting the singular source terms, discrete Green functions on finite lattices with homogeneous and Neumann boundary conditions,  a new discovered simplest $L^2$ discrete delta function,  and the convergence proof of the IB method for
elliptic interface problems \cite{li:mathcom}.
The conclusion in this paper can apply to problems with different  boundary conditions  as long as the problems are wellposed.
The proof process
also provides an efficient way to decouple the system into three Helmholtz/Poisson equations without affecting the order of  convergence.
A non-trivial numerical example is also provided to confirm the theoretical analysis and the simple new discrete delta function.
\end{abstract}



\section{Introduction}

 The Immersed Boundary (IB) method \cite{peskin:heart2}  is one of the most important modeling tools and numerical methods. The IB method has been applied  to  many   problems in
  mathematics and  engineering including  biology, fluid mechanics, material science, electric-magnetics, magnetohydrodynamics,
  see for example, \cite{Mittal_IB,PM95} for  reviews and references therein.

\begin{figure}[hpbt]
 \centerline{  \includegraphics[width=0.40\textwidth]{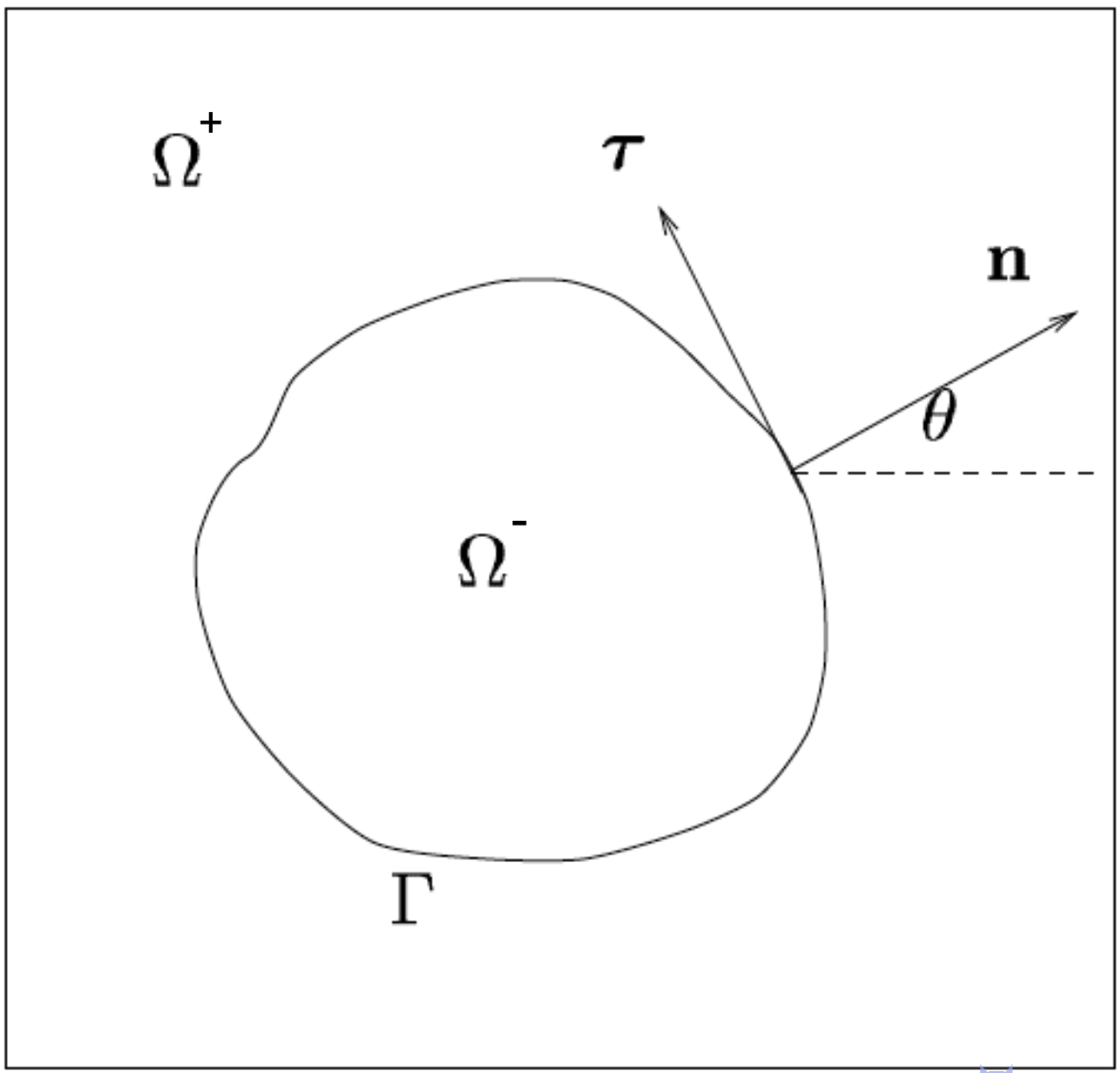} }
 \caption{A diagram of the Peskin's Stokes   IB model with an immersed  interface
 $\Gamma$.} \label{domainD}
\end{figure}

 Consider the stationary incompressible Stokes equations with an immersed interface $\Gamma$ as in the Peskin's IB model and a prescribed velocity boundary condition,
 \eqm
 && \nabla p = \mu \Delta  {\bf u}  + {\bf G} + \int_{\Gamma} {\bf f}(s)\,
\delta_2({\bf x}-{\bf X}(s)) ds, \quad {\bf x} \in \Omega, \label{Sa}\\ \eqsp
&& \nabla \cdot {\bf u} = {0}, \quad {\bf x} \in \Omega, \label{Sb} \\ \eqsp
&& \left. {\bf u} ({\bf x})\right |_{\partial \Omega}={\bf u} _{0}({\bf x}). \label{Sc}
\enm
This paper will focus on two-dimensional (2D) problems. Thus, in the expression above, we have ${\bf x}=(x,y)$, ${\bf u} =(u,v)$,   ${\bf f}=(f_1,f_2)$, and so on.
 We assume that  $\mu$ is a constant,   the external bounded source ${\bf G} \in C(\Omega^{\pm})$,  the source strength  $ {\bf f}(s)\in C^1(\Gamma)$,  and the interface $\Gamma(s)={\bf X}(s)=(X(s),Y(s))$ with a parameter $s$, {\em e.g.}, the arc-length, see Fig.~\ref{domainD} for an illustration.
 In the error analysis, the error of the velocity at the boundary is zero. Thus, for simplicity, we can simply use the term 'no-slip' boundary condition.  

Alternatively, the incompressible Stokes equations can be written as two regular Stokes equations on two separated domain but couple with jump conditions, see \cite{rjl-li:sinum,li:thesis,li:book}, as 
\eqm  \label{SaB}
 && \nabla p = \mu \Delta  {\bf u}  + {\bf G},  \quad {\bf x} \in \Omega\setminus \Gamma, \\ \eqsp
&& \nabla \cdot {\bf u} = 0, \quad {\bf x} \in \Omega, \label{Sb} \\ \eqsp
&& \left. {\bf u} (x,y)\right |_{\partial \Omega}={\bf u} _{0}(x,y), \\ \eqsp
&& [p] =  \hat{f}_1(s), \quad [p_n] = \frac{\partial }{\partial  s}  \hat{f}_2(s) + [{\bf G} \cdot  {\bf n}] , \\ \eqsp
 && [ {\bf u} ] = {\bf 0} , \quad \left [  {\bf u}_ n  \right ] =  \hat{f}_2(s)  {\bftau},  \label{SeB}
\enm
where ${\bf G}=[G^1, \, G^2]^T$, ${\bf n}$ and ${\bftau}$ are the unit normal and tangential directions, respectively,  $[p_n] = [\nabla p \cdot {\bf n} ]=[\frac{\partial p}{\partial  n}]$, for example, is the normal derivative of the pressure and so on. 
The jumps, for example, $[p_n]({\bf X})$, ${\bf X}\in \Gamma$,  is defined as
\eqm
  [p_n]({\bf X} ) = \lim_{{\bf x}\goto {\bf X}, {\bf x}\in \Omega^+} \frac{\partial p}{\partial n} ({\bf x}) 
   - \lim_{{\bf x}\goto {\bf X}, {\bf x}\in \Omega^-} \frac{\partial p}{\partial n} ({\bf x}),
\enm
and so on. 
The normal and tangential  force densities are defined as
\eqm
  \hat{f}_1(s) = {\bf f}(s) \cdot {\bf n}, \qquad \hat{f}_2(s)=  {\bf f}(s) \cdot {\bftau}.
\enm
In component-wise form, we also have
\eqml{un_p}
   [u_n] &=& [p]\cos \theta - f_1 =\hat{f}_2 \sin \theta,\\ \eqsp
   [v_n] &=& [p]\sin \theta - f_2= -\hat{f}_2 \cos \theta ,
\enml
where $\theta$ is the angle between the normal direction and the $x$-axis, {\em cf.} Figure~\ref{domainD}.
Advantages of the expressions above include {\em getting rid of singular Dirac delta functions} in terms of the jump conditions, also called internal boundary conditions so that high order methods such as IIM, IFEM, MIB,  the virtual nodes method, and other high order methods can be applied. And more importantly, the different singularities for the primitive variables and relations  to  the source strength, which is important in our convergence proof, can be revealed in the jump conditions. With the above representation,  the local truncation errors of the IB method can {\em be treated as a regular problem} without  Dirac delta functions similar to the role of a Peskin's discrete delta function, which approximates a singular source term with  a `regular' function.

 We assume that ${\bf u}\in C^2(\Omega^{\pm})\cap C(\Omega) $ and $p({\bf x} ) \in C^1(\Omega^{\pm})$ throughout the paper, that is, the velocity is continuous in the entire domain,  and the first and second partial derivatives are piecewise continuous and bounded excluding the interface; the pressure and its first partial derivatives are piecewise continuous and bounded excluding the interface.

While there are various implementations,   there are two prominent  IB methods for solving the Stokes equations with a prescribed velocity boundary condition. One is based on a MAC scheme in which the discrete pressure and velocity are coupled together. Advantages of the MAC approach include no need to provide numerical boundary conditions for the pressure,  and that the computed velocity is divergence free in a discrete sense. The challenge is how to solve the resulting linear system of equations,  which is a saddle problem,  efficiently.  The other popular IB method is the so-called three Poisson equations approach in which the pressure is solved first with an approximate Neumann boundary condition. The velocity then can be obtained by solving other two Poisson equations. The most important advantage is that fast Poisson solvers can be utilized. The discrete divergence of the computed velocity
is usually of $O(h^2)$.  In this paper, the two methods are closely related as we can see later.


It is known that Peskin's IB method is  first order accurate for the velocity  and zeroth order for the pressure in the pointwise (infinity) norm.
However, there was almost no rigorous proof  in the literature until
\cite{mori-proof}, in which the author has proved the first order accuracy of the velocity of the IB method for the Stokes equations with a periodic boundary condition. The proof is based on some known inequalities between the fundamental solution and the discrete Green function with a periodic boundary condition for Stokes equations. In \cite{lai_IB_Proof},   the author pointed out that
the pressure obtained from the IB method has $O(h^{1/2})$ order of convergence in the $L^{2}$ norm without a proof.

In \cite{li:mathcom}, the first order convergence of the IB method was  proved for elliptic interface problems with a Dirichlet boundary condition using the discrete Green function on bounded lattices. The main goal of this paper is to provide a convergence proof for the IB method for Stokes equations with a prescribed velocity along the boundary.  One of  the difficulties in the proof is that the pressure and the velocity are coupled together.
We discuss the proof for two commonly used IB methods  for Stokes equations.  
The first IB method is the three Poisson equations approach, see \cite{rjl-li:stokes,xwan-li-pbc,li-ito-lai,xwan-li-level} and Section~\ref{sec:3poisson}. The second IB method the classical  one using a MAC grid, see for example \cite{IB_Griffith_07,Bringley08,KOLAHDOUZ_2020,Stein-Guy-19} and the references therein.
Our strategy is to get an error estimate for the pressure first. We show that even though the pressure has $O(1)$ error near the interface, it is $O(\sqrt{h |\log h |})$  convergent in the $L^2$ norm. Furthermore,  the pressure has $O(h)$  pointwise convergence $\sqrt{h}$ away from the interface.  Then, we show that the singularity from the discontinuous pressure is offset by  a part of singular source along the interface, and thus the proof process for the scalar elliptic interface problem can be applied.

Note that, many dedicated and fine analysis about the convergence of numerical methods using a MAC grid for incompressible Stokes and Navier-Stokes equations can be found in the literature, see for example \cite{MR3342722,MR1829554,Nicolaides-92,MR3648073}. Nevertheless, the focus of this paper is about the convergence analysis when the  IB methods are applied to Stokes equations with  a prescribed velocity boundary condition for which many analysis intended for high order or super-linear convergence do not apply or are not needed.

The rest of the paper is organized as follows. In the next section, we introduce discrete Green functions on finite lattices with homogeneous Dirichlet and Neumann boundary conditions. In Section~\ref{sec:3poisson}, we introduce the three Poisson equations for solving the Stokes equations and the convergence analysis. A new simplest discrete delta function in $L^2$ is introduced as needed in the proof. In Section~\ref{sec:MAC}, we introduce the IB method for solving the Stokes   equations using mark-and-cell (MAC) meshes and the convergence analysis.   In Section~\ref{sec-ex}, we use a nontrivial example in which both  the pressure and the flux of the velocity have non-constant jumps across a circular interface to confirm the theoretical analysis.

\section{Discrete Green functions on bounded grids}

Discrete Green functions play  important roles in the convergence proof.  The discrete Green function on a infinite lattatice centered at one particular grid has been well studied in the literature, see for example \cite{rutka-phd-thesis,thomee} and the reference therein.  However, there are limited discussions on discrete Green functions on bounded domains. A early work on construction and estimates can be found in \cite{rutka-phd-thesis} in which the discrete Green function  with a homogeneous Dirichlet boundary conditions is defined below.


\subsection*{Dirichlet discrete Green function on a bounded grid}

\begin{definition}
Given a uniform grid, ${\bf x}_{ij}$, $i,\, j = 0, 1, \cdots, N$. 
 Let ${\bf e}_{lm}$, $1\le i \le N-1$,  $1\le j \le N-1$,  be the unit grid function whose values are zero at all grid points except
  at ${\bf x}_{lm}=(x_l, y_m)$ where its component is $e_{lm}=1$.  The discrete Green function centered at ${\bf x}_{lm}$ with
   homogeneous boundary condition   is defined as
  \eqm \label{discreG}
     {\bf G}^{h} \lp {\bf x}_{ij}, {\bf x}_{lm}\rp = \lp A_h^{-1} {\bf e}_{lm} \frac{1}{h^2} \rp_{ij}, \qquad  {\bf G}^{h} \lp \partial \Omega_h, {\bf x}_{lm}\rp=0,
  \enm
  where $A_h$ is the standard centered five-point discrete Laplacian, 
\eqml{discreG1}
  &   \dsp A_h U_{ij} =   \frac{ U_{i-1,j} + U_{i+1,j} + U_{i,j-1}+U_{i,j+1}-4 U_{i,j}}{h^2}=\left \{ \begin{array}{ll}
  \dsp \frac{1}{h^2} &  ~ \mbox{if~  $i=l$, $j=m$, } \\ \eqsp
  \dsp  0  &    \mbox{otherwise,}
 \end{array} \right. \\ \eqsp
 & \dsp U_{0j}=U_{Nj}=U_{i0}=U_{iN}=0, \qquad i,j=0,1,\cdots, N.
\enml
We use $\partial \Omega_h$ to denote all the boundary grid points.
\end{definition}

The discrete Green function defined above is unique since $-A_h$ is a symmetric definite matrix. It has ben shown in \cite{} that ${\bf G}^{h}$ can be expressed as a  finite linear combination of  the Green functions ${\bf g}^{h} \lp {\bf x}_{ij}, {\bf x}_{l_1,m_1}\rp $ on an infinite lattice. ${\bf G}^{h} \lp {\bf x}_{ij}, {\bf x}_{lm}\rp $ behaviors just like the ${\bf g}^{h} \lp {\bf x}_{ij}, {\bf x}_{lm}\rp$ with similar estimates in terms of the grid function and its partial derivatives if the center is sufficiently away from the boundary. 

Using a high order interpolation, say in $C^1(\Omega)\cap H^2(\Omega)$, we can extend the discrete Green function  ${\bf G}^{h} \lp {\bf x}_{ij}, {\bf x}_{lm}\rp$ to the entire domain as ${\bf G}^{h} \lp {\bf x}, {\bf x}_{lm}\rp$. 
In general, there are infinite interpolation functions. 
We consider any such an interpolation function that satisfies the following, see \cite{} for a constructing process,
\label{ghi_cond}
\begin{itemize}
   \item ${\bf G}^{h}_I({\bf x}_{ij},{\bf x}_{lm})={\bf G}^{h}({\bf x}_{ij},{\bf x}_{lm})$.
  \item ${\bf G}^{h}_I({\bf x},{\bf x}_{lm})\in C^1(\Omega)\cap H^2(\Omega)$.

  \item $\Delta_h {\bf G}^{h}_I({\bf x}_{ij},{\bf x}_{lm})=\Delta_{h}{\bf G}^{h}({\bf x}_{ij},{\bf x}_{lm})=0$, that is, zero for all $i$ and $j$ except for $i=l$ and $j=m$.


  \item The interpolation is third or higher order accurate, that is
  \begin{equation}\label{order}
    \left |\partial^{\mathbf{\alpha}} {\bf G}^{h} \left ({\bf x}_{ij},{\bf x}_{lm} \right) - \partial^{\mathbf{\alpha}} {\bf G}^{h}_I \left ({\bf x}_{ij},{\bf x}_{lm} \right) \right |\le
    C h^{3-\|\mathbf{\bfalf}\|_1}, \quad \mbox{for} \quad \|\mathbf{\bfalf}\|_1\le 2,
  \end{equation}
  where the derivatives of ${\bf G}^{h}  ({\bf x}_{ij},{\bf x}_{lm})$ is defined from finite differences, see \cite{thomee}, and ${\mathbf{\alpha}}$ is the multi-index notation as used in the literature for the Sobolev spaces.

   \item  $\dsp \iint_{R_{lm}}\Delta {\bf G}^{h}_I({\bf x},{\bf x}_{lm})dx dy=O(h)$ except for the four neighboring squares centered at ${\bf x}_{lm}$ on which $\dsp \iint_{R_{lm}}\Delta {\bf G}^{h}_I({\bf x},{\bf x}_{lm})dxdy=1 + O(h)$, where $R_{lm}$ is the square $[x_l-h,\; x_l+h]\times [y_m-h,\; y_m+h]$ and so on.
\end{itemize}

\vthin

With the above interpolation interpolation Green functions, the following estimates hold.
\begin{lemma}
 Let ${\bf G}^{h}_I ({\bf x},{\bf x}_{lm} )$ be an interpolation function of
 ${\bf G}^{h}({\bf x}_{ij},{\bf x}_{lm} )$ that satisfies the conditions above, then we have
the following estimates:
 \begin{eqnarray}
   &&\dsp \frac{}{}{\bf G}^{h}_I \left ({\bf x},{\bf x}_{lm} \right) \le \frac{1}{4} +
    \frac{1}{16} \log \left ( \| {\bf x}-{\bf x}_{lm}\|_2^2 + h^2 \right) + O(h), \label{u_bd}\\ \eqsp
   && \left| \partial^{\mathbf{\alpha}} {\bf G}^{h}_I \left ({\bf x},{\bf x}_{lm} \right)\right |\le
     \frac{C}{ \dsp \frac{}{}\left( \| {\bf x}-{\bf x}_{lm} \|_2  + h\right)^{\|\mathbf{\bfalf}\|_1}} + O(h),
     \quad \mbox{if $\|\mathbf{\bfalf}\|_1\le k-1$},
     \label{deri_est}
 \end{eqnarray}
 where $k$ is the order of interpolation and $C$ is a constant. The inequality \eqn{deri_est} is true if $\;dist({\bf x}, \partial \Omega)\sim O(1)$ and $\;dist({\partial \Omega}, {\bf x}_{lm}) \sim O(1)$.
\end{lemma}

The estimates above are expected since the discrete Green function is an approximation to the green function $\Delta G = \delta ({\bf x} - {\bf x}_{lm})$ when ${\bf x}$ is away from the center and the boundary. 

\begin{figure}[hpt]
\begin{minipage}[t]{3.0in} (a)

\centerline{\includegraphics[width=0.85\textwidth]{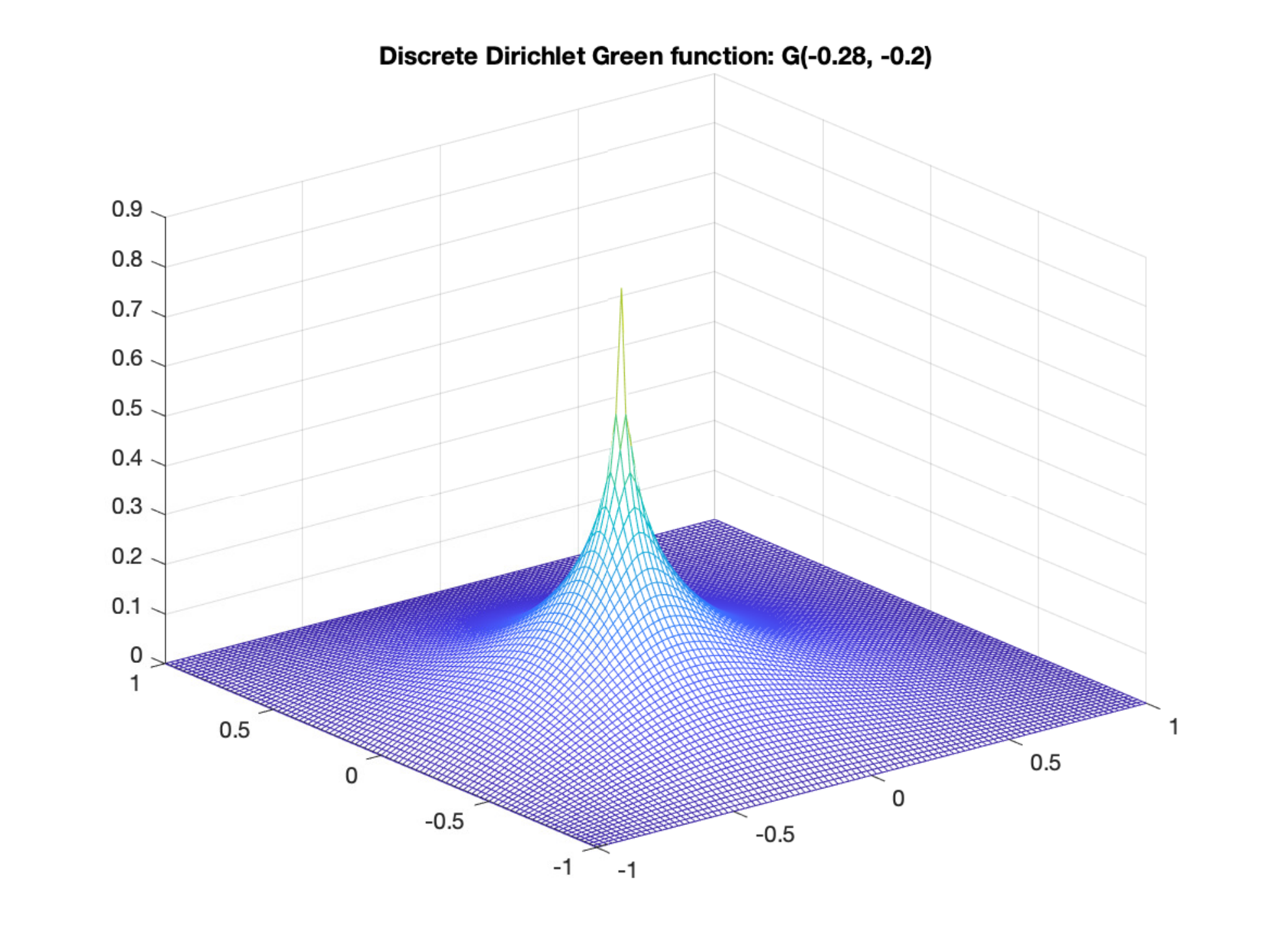}}
 \end{minipage}
\begin{minipage}[t]{3.0in} (b)

\centerline{\includegraphics[width=0.85\textwidth]{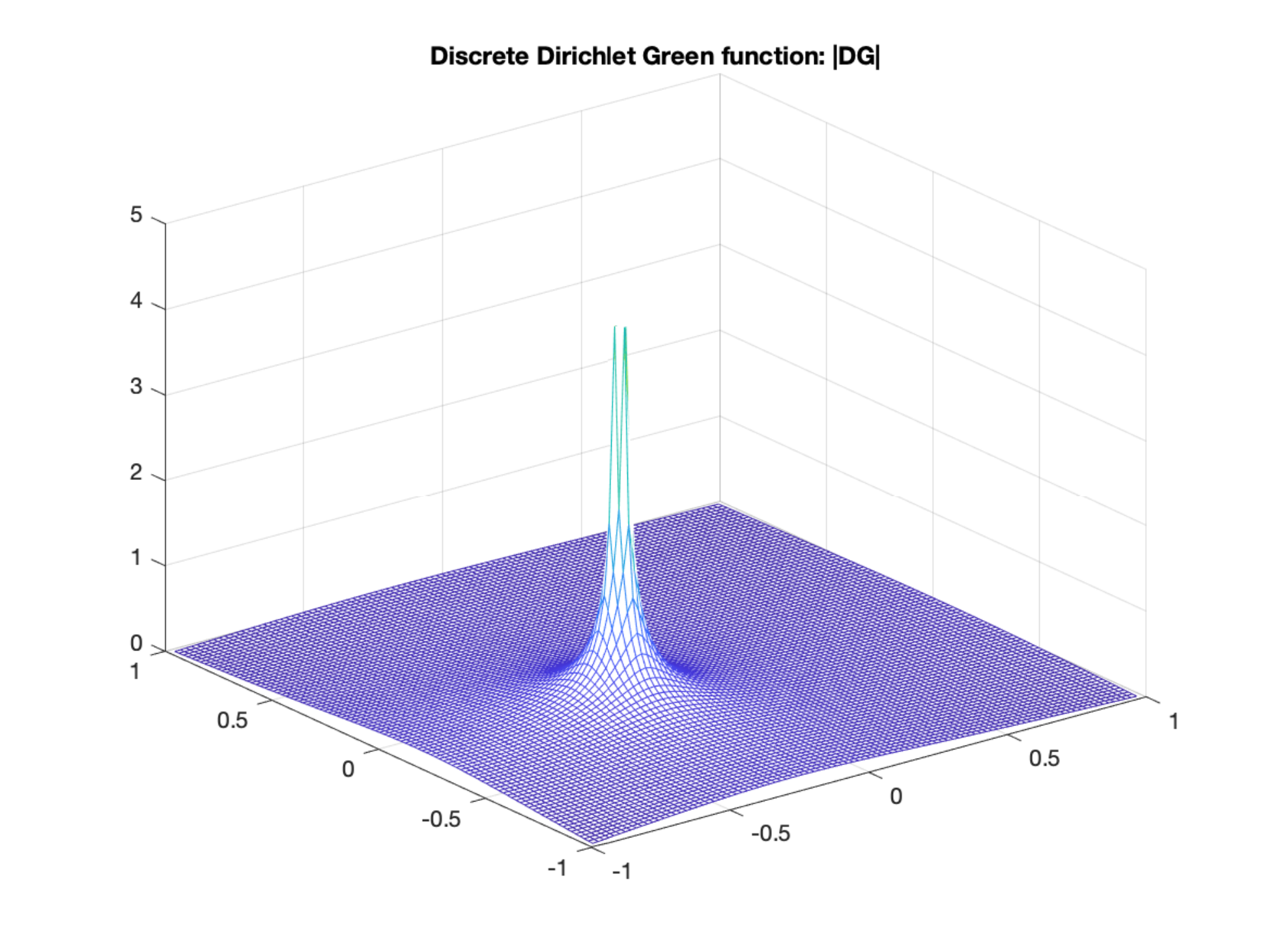}}
 \end{minipage}

\begin{minipage}[t]{3.0in} (c)

\centerline{\includegraphics[width=0.85\textwidth]{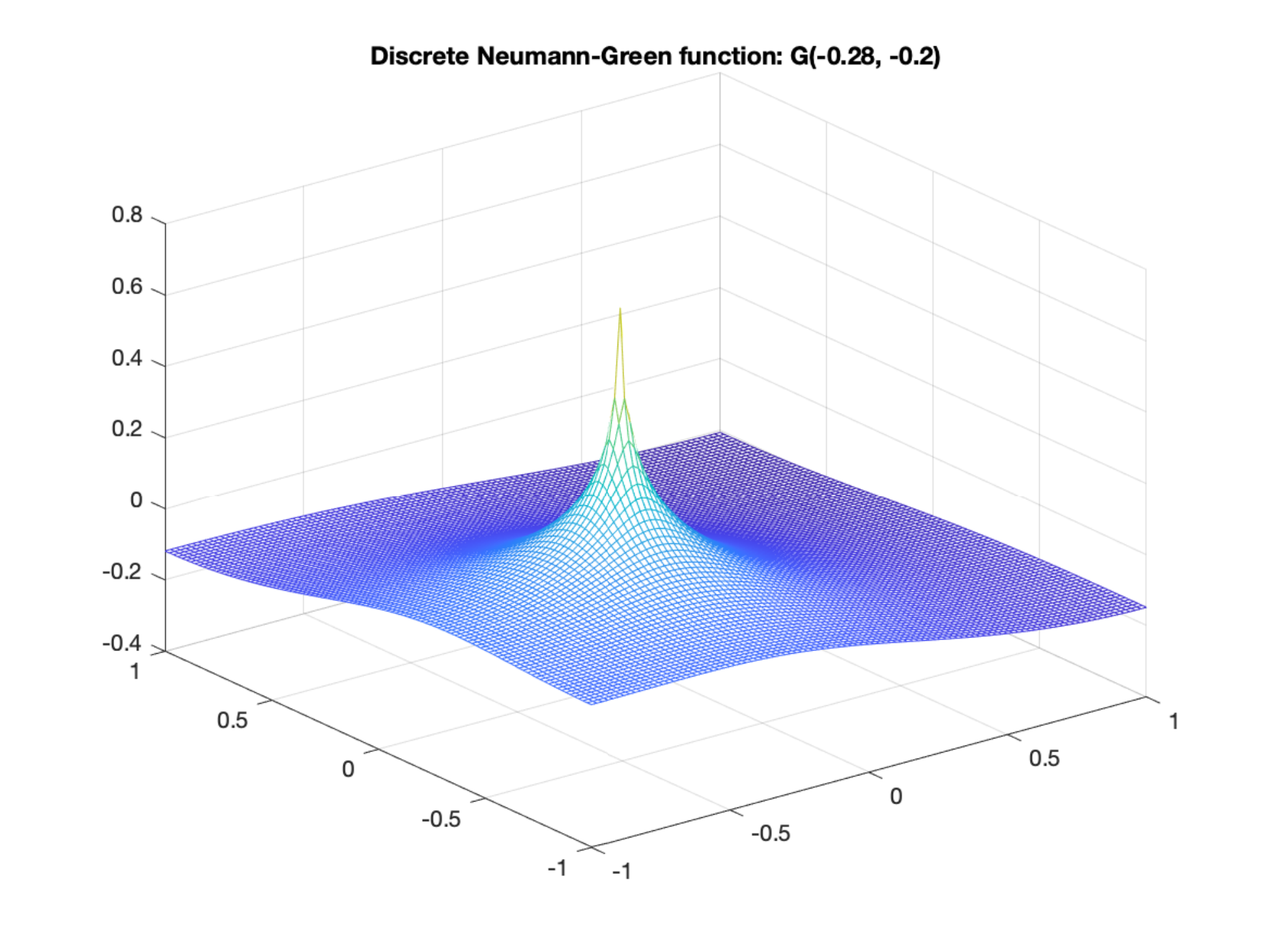}}
 \end{minipage}
\begin{minipage}[t]{3.0in} (d)

\centerline{\includegraphics[width=0.85\textwidth]{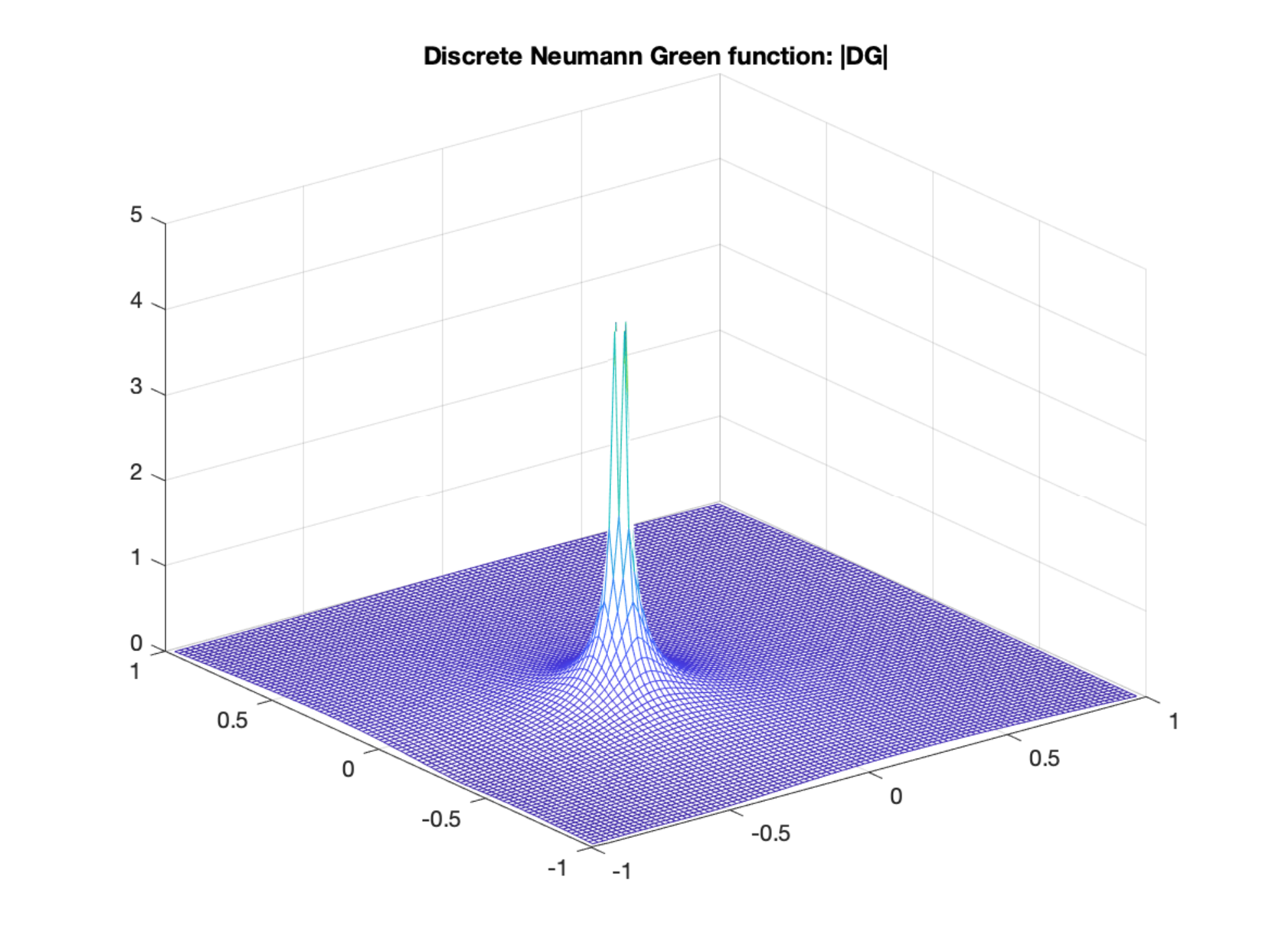}}
 \end{minipage}

\caption{Discrete Green functions centered at ${\bf x}_{lm} = (-0.28,-0.2)$  and their magnitude of the first order discrete derivatives.  Dirichlet BC, (a)-(b):  $ - {\bf {G}}^{h} \lp {\bf x}_{ij}, {\bf x}_{lm}\rp$ and $ \left | \nabla_h {\bf {G}}^{h} \lp {\bf x}_{ij}, {\bf x}_{lm}\rp \right |$;   Neumann BC,
(c)-(d):  $ - {\bf \hat{G}}^{h} \lp {\bf x}_{ij}, {\bf x}_{lm}\rp$ and $ \left |  \nabla_h {\bf \hat{G}}^{h}  \lp {\bf x}_{ij}, {\bf x}_{lm}\rp \right |$, where $\nabla_h$ is the discrete gradient using the central finite difference formula. \label{fig:Green}
}
\end{figure}

\subsection*{Nuemann discrete Green function on a bounded grid}

Similar to the Dirichlet discrete Green function,  the Neumann discrete Green function centered at ${\bf x}_{lm}$ with
   homogeneous flux boundary conditions ide defined below.
\begin{definition}
 Given a uniform grid, ${\bf x}_{ij}$, $i,\, j = 0, 1, \cdots, N$. 
 Let ${\bf e}_{lm}$ be the unit grid function whose values are zero at all grid points except
  at ${\bf x}_{lm}=(x_l, y_m)$ where its component is $e_{lm}=1$.  The discrete Neumann Green function centered at ${\bf x}_{lm}$ with
   homogeneous boundary condition   is defined as
\eqml{discreG}
  &   {\bf \hat{G}}^{h} \lp {\bf x}_{ij}, {\bf x}_{lm}\rp = \lp \hat{A}_h^{-1} {\bf e}_{lm} \frac{1}{h^2} \rp_{ij}, \quad i,j=0,1,\cdots,N,
\enml
where the Neumann discrete Laplacian $\hat{A}_h$ is defined as
\eqml{discreG2}
  &   \dsp \hat{A}_h U_{ij} =   \frac{ U_{i-1,j} + U_{i+1,j} + U_{i,j-1}+U_{i,j+1}-4 U_{i,j}}{h^2}=\left \{ \begin{array}{ll}
  \dsp \frac{1}{h^2} &  ~ \mbox{if~  $i=l$, $j=m$, } \\ \eqsp
  \dsp  0  &    \mbox{otherwise,}
 \end{array} \right.
\enml
in the interior, $i,\ j = 2,\cdots, N-2$, and  the
discrete homogeneous Neumann boundary condition is defined at the boundary grid points,
\eqml{discreG3}
    & \dsp  \frac{ U_{2,j}  - U_{1,j} }{h} = 0,  \qquad  \frac{ U_{N,j}  - U_{N-1,j} }{h} = 0,  \qquad j=1,\cdots, N, \\ \eqsp
    &  \dsp  \frac{ U_{i,2}  - U_{i,1} }{h} = 0,  \qquad  \frac{ U_{i,N}  - U_{i,N-1} }{h} = 0,  \qquad i=1,\cdots,  N.
\enml
To make the solution to the above linear system of equations unique, we also set
$ U_{11} = 0$.
\end{definition}

The discussion of the discrete Neumann Green function  on a bounded grid can be carried out almost in the same way as that of the Dirichlet Green function on a bounded grid. 
In Fig.~\ref{fig:Green},  we plot the two different (Dirichlet and Neumann) Green functions and the magnitude of  their discrete gradient using the standard centered finite difference  formula. We see clearly similar behaviors especially near the center of the point source.


\section{Convergence analysis of IB method of the three Poisson equations approach} \label{sec:3poisson}.

The simplest IB method for solving the incompressible Stokes flows with a constant viscosity is to use the  three Poisson equations approach on uniform meshes, see for example \cite{rjl-li:stokes}. 
Applying the divergence operator to the momentum equations \eqn{Sa}, we obtain
\eqm
  \Delta p = \grad \cdot {\bf G} +  \int_{\Gamma} \grad \cdot \lp \frac{\null}{\null} {\bf f}(s)\,  
\delta_2({\bf x}-{\bf X}(s)) \rp ds. 
\enm
An approximate boundary condition is needed since the pressure boundary condition is not a free variable.  A common practice is to use the homogeneous Neumann boundary condition $\left . \frac{\partial p}{\partial n} \right |_{\partial \Omega}= 0$. The pressure can differ by a constant. Usually we set a reference pressure, say $p(a,c)=0$ to make it unique. Once the pressure has been solved, we can solve the two Poisson equations from the momentum equation to obtain the velocity. 

Given a uniform grid,
\eqm
 x_i = a +  i h, \quad i=0,1,\cdots, N;  \qquad y_j= c   + j h , \quad j=0,1,\cdots, N,
\enm 
along with a number interface points, ${\bf X}_k$, $k=1,2,\cdots, N_b$ with the assumption that $\dsp \max_k s_k = \max_k \|{\bf X}_{k+1}-{\bf X}_{k}\|_2 \sim O(h)$. 
The approximate solution to the pressure is the solution to the following system of equations
\eqml{plap}
 &  \dsp  \frac{P_{i-1,j}+P_{i+1,j}+P_{i,j-1}+P_{i,j-1 }-4P_{i,j }}{h^2}  = D_{ij}^x \lp G^1_{ij} + F^1_{ij} \rp + D_{ij}^y \lp G^2_{ij} + F^2_{ij} \rp, 
\enml
where $i, j=1, \cdots, N-1$, $G^1_{ij}$ and $D_{ij}^y$ are the central finite operators defined as
\eqm  
 D_{ij}^x  G^1_{ij} = \frac{G^{1}_{i+1,j} - G^{1}_{i-1,j}}{2 h}, \qquad D_{ij}^y G^2_{ij}=   \frac{G^{2}_{i,j+1} - G^{2}_{i,j-1}}{2 h}, 
\enm
and $F^1_{ij}$ and $F^2_{ij}$ are the source distribution from a consistent discrete delta function, for example,
\eqm
 F^1_{ij} = \sum_{k=1}^{N_b}  f_1({\bf X}_k) \,  \delta_h(x_i- X_k) \, \delta _h(y_j- Y_k) \Delta s_k, \quad \Delta s_k = \left \| \frac{\null}{\null} {\bf X}_{k+1} - {\bf X}_k \, \right  \|_2, 
\enm
and so on. 
The homogeneous Neumann boundary condition can be simply discretized with the  forward finite difference method. 
\eqmno 
  \frac{P_{1,j}-P_{0,j}}{h}=0, \quad  \frac{P_{N,j}-P_{N-1,j}}{h}=0, \quad  \frac{P_{i,1}-P_{i,0}}{h}=0, \quad  \frac{P_{i,N}-P_{i,N-1}}{h}=0.
\enmno

Once we have computed the pressure, we can obtain the velocity for the momentum equation, which are two Poisson solvers for the velocity components,
\eqml{plap-unif}
 &  \dsp  \mu \frac{U_{i-1,j}+U_{i+1,j}+U_{i,j-1}+U_{i,j+1 }-4U_{i,j }}{h^2}  = \frac{P_{i+1,j}- P_{i-1,j}}{2 h} - F^{1}_{i,j} ,  \\ \eqsp
 & i=1, \cdots, N-1, \qquad j=2,\cdots,  N-1,
\enml
and
\eqml{plap-unif}
 &  \dsp  \mu \frac{V_{i-1,j}+V_{i+1,j}+V_{i,j-1}+V_{i,j+1 }-4V_{i,j }}{h^2}  = \frac{P_{i,j+1}- P_{i,j-1}}{2 h} - F^{2}_{i,j} ,  \\ \eqsp
 & i=1, \cdots, N-1, \qquad j=2,\cdots, ... , N-1,
\enml
with the prescribed Dirichlet boundary conditions. 

There are several advantages of the three Poisson equations approach. The pressure and velocity are decoupled; fast Poisson solvers  based on uniform meshes can be utilized so the method is fast; The implementation of the IB method using a discrete delta function will yield first order accurate velocity, which is comparable to the IB method using a MAC grid and will be proved in this paper. The disadvantage is that the computed velocity is approximately divergence free. 

\subsection{Convergence analysis of  the pressure from the three Poisson equation approach}

In this paper, we will use a generic error constant $C$ in all  the estimates. Without loss of generality, we assume that $\mu=1$, 
${\bf u}|_{\partial \Omega}={\bf 0}$, and ${\bf G}={\bf 0}$.
We define errors in pressure, velocity as
\eqm
  E^u_{ij} = u(x_i,y_j) - U_{ij}, \qquad E^v_{ij} = v(x_i,y_j) - V_{ij}, \qquad E^p_{ij} = p(x_i,y_j) - P_{ij}.
\enm

It is known that the pressure computed from IB methods does not have point-wise convergence in general due to the dipole singularity of the source term. However, we prove $O(h\log h)$ convergence in the $L^2$ norm, and pointwise $O(h\log h)$ convergence  for grid points that are $O(\sqrt{h})$ away from the interface. 

We call a grid point $(x_i,y_j)$ as a {\em regular} one for the pressure  discretization if the interface does not cut through the standard centered five-point stencil, or in other words, all the five grid points are from the same side of the interface. Otherwise the grid point is called {\em irregular}.

Define the local truncation errors of the pressure from the three-Poisson equation approach as 
\eqml{trunp}
 &  \dsp  T_{ij}^p = \frac{p(x_{i-1},y_j) +p(x_{i+1},y_j)  +p(x_i,y_{j-1}) +p(x_i,y_{j+1}) -4p(x_i,y_j)}{h^2}  - D_{ij}^x   F^1_{ij}  - D_{ij}^y F^2_{ij} , 
\enml
where $i, j=1, \cdots, N-1$, and at the boundary
\eqml{ptrun-bc}
  \dsp T_{0,j}^p = \Delta p(x_0,y_j) -  \frac{p(x_1,y_j) -p(x_0,y_j) }{h},  && \dsp T_{N,j}^p =  \Delta p(x_N,y_j) -  \frac{p(x_N,y_j) -p(x_{N-1},y_j) }{h},  \\ \eqsp
   \dsp   T_{i,0}^p =  \Delta p(x_i,y_0)  - \frac{p(x_i,y_1) -p(x_i,y_0) }{h},   &&  \dsp T_{i,N}^p =  \Delta p (x_i,y_N)  - \frac{p(x_i,y_N) -p(x_i,y_{N-1}) }{h}.
\enml

It is easy to check that the local truncation errors of the pressure  are $O(h^2)$ at regular interior grid points;  $O(1/h^2)$ at irregular grid points due to the jump condition in the pressure, and $O(1)$ at the boundary points, which is shown below.
Take  $(x_0,y_j)$,  for example,  the local truncation error for the pressure  satisfies
 \eqmno
  T_{0,j}^p   &=& \Delta p (x_1,y_j) - \lp \frac{p(x_0, y_j )  - p(x_1,y_j ) }{h}  \rp  \\ \eqsp
   &\sim &   O(1) +  p_x( x_0, y_j )  + O(h)  \sim  O(1).
 \enmno
 In the above, we have used  the fact that the pressure is smooth away from the interface up to the boundary, and the forward  finite difference approximation is first order accurate. Now we are ready to prove the convergence of the pressure. 

Let $W$ be the maximum number of layers of the discrete delta function, that is $\delta_h({\bf x} )=0$ if 
 $| {\bf x} | \ge W h$. For example, $W=1$ for the hat discrete delta function, while $W=2$ for discrete cosine discrete delta function.
  From  \eqn{trunp}-\eqn{ptrun-bc}, we can write down the error of the pressure as a matrix-vector form,
 \eqm
  \hat{A}_h \, {\bf E}^p  = {\bf T}^{p},
 \enm
 where ${\bf T}^{p}$ is composed of the local truncation errors that are $O(h^2)$ at regular interior grid points;  $O(1/h^2)$ at irregular grid points due to the jump condition in the pressure, and $O(1)$ at the `boundary points'.  Thus, we derive,
 \eqml{eplm}
 \hat  E^p_{lm} &=& \dsp \lp (\hat{A}_h)^{-1} {\bf T}^{p} \rp_{lm}   \\ \eqsp
 &=&  \lp (\hat{A}_h)^{-1} {\bf T}_{reg}^{p} \rp_{lm} + \lp (\hat{A}_h)^{-1} {\bf T}_{irr}^{p} \rp_{lm} + \lp (\hat{A}_h)^{-1} {\bf T}_{BC}^{p} \rp_{lm} \\ \eqsp
  &=&  \dsp O(h^2) + \lp (\hat{A}_h)^{-1}  {\bf T}_{irr}^{p} \rp_{lm} + O(h)  \\  \eqsp
  &=& \dsp \sum_{ dist({\bf x}_{ij}, \Gamma)\le \hat{W} h
  }\lp h^2 T_{ij}   (\hat{A}_h)^{-1} {\bf e}_{ij} \frac{1}{h^2} \rp_{lm} + O(h) + O(h^2) \\ \eqsp
  &=&  \dsp \sum_{ dist({\bf x}_{ij}, \Gamma)\le \hat{W} h }  h^2 \lp \Delta_{h} p(x_i,y_j) - \hat{C}_{ij}^{IB}\rp {\bf \hat{G}}^{h}({\bf x}_{ij},{\bf x}_{lm}) + O(h) + O(h^2) \\ \eqsp
  &=& \dsp \sum_{ij}  h^2  \Delta_{h} p(x_i,y_j) {\bf \hat{G}}^{h}({\bf x}_{ij},{\bf x}_{lm}) - \sum_{ij}  h^2 \hat{C}_{ij}^{IB}  {\bf \hat{G} }^{h}({\bf x}_{ij},{\bf x}_{lm})  + O(h)+ O(h^2), \\ \eqsp
   &=& \dsp \lp \int_{\Gamma} (\grad \cdot {\bf f}(s))  \, {\bf \hat{G}}^h_I({\bf X}(s),{\bf x}_{lm}) ds - \sum_k   (\grad \cdot {\bf f}(s_k)) {\bf \hat{G}}_I^h ({\bf X}_k),{\bf x}_{lm}) \Delta s_k\rp + O(h) . 
\enml
In the expressions above
$ {\bf \hat{G}}^h({\bf x}, {\bf x}_{lm})$ is the  discrete Neumann Laplacian centered at ${\bf x}_{lm})$, $\bar{C}_{ij}^{IB}$ is the result of the discrete delta function that distributes the singular source term to  the grid points near the interface~$\Gamma$. We have used the consistence of the discrete delta function in the interpolation and source distribution to obtain the last expression. We have $\hat{W}=W+2$ due to the divergence operator applying to the source term. 

The above expression of $E^p_{lm}$ is almost the same as that in \cite{li:mathcom} for a Poisson equation with a closed interface except that the source here is a dipole, or $O(1/h)$ more singular. Thus, the estimates there are still valid with an  $O(1/h)$ factor, some of which are summarized in the following lemma. 

\begin{lemma}
Assume that $p(x,y)\in C^2(\Omega^{\pm})$ and $h$ is small enough so that $dist(\Gamma,\partial \Omega)\sim O(1)$ and the area enclosed by $A_{\Omega^-} \sim O(1)$, then the computed pressure from the three Poisson equation approach satisfies the following error estimate.
\bi
\item If \ $l,m \neq 1$ or $n-1$, then we have,
\eqm
\sum_{ij}  \Delta_{h} p(x_i,y_j)  \hat{{\bf G}}^{h}({\bf x}_{ij},{\bf x}_{lm}) h^2 =
 \int_{\Gamma}  (\grad \cdot {\bf f}(s))   \hat{{\bf G}}^h_I ({\bf X}(s),{\bf x}_{lm}) ds + O(1).
\enm
\item Let $\hat{C_{ij}}^{IB}$ be the  non-zero contributions from the distribution of \ $\grad \cdot {\bf f}(s)$ in the immersed boundary method, then 
\eqmno
\sum_{ij}  \hat{C_{ij}}^{IB}  \hat{{\bf G}}^{h}({\bf x}_{ij},{\bf x}_{lm}) h^2 =
 \int_{\Gamma}   (\grad \cdot {\bf f}(s))   \, \hat{{\bf G}}_I^h \!\!\left({\bf X}(s),{\bf x}_{lm} \frac{\null}{\null}\!\!\right) \, ds + O( \log h).
\enmno
\item  If \ $dis({\bf x}_{lm},\Gamma) \sim O(1)$, then
\begin{eqnarray*}
 \left  |E^p_{lm} \right |   & \leq &   \frac{C h}{\max_{k} \lp \|{\bf X}_k- {\bf x}_{lm}\|_2 + h\rp^2 \frac{}{}} \sim \frac{ C h}{\lp O(1) + h\rp^2} \leq Ch,
\end{eqnarray*}
\item If \ $dist({\bf x}_{lm}, \Gamma)\ge  \bar{W} h$, where $\bar{W}= W+2$,  then 
\eqm
 \left  | E^p_{lm} \right |  \le  \frac{C h^2}{ (dist( {\bf x}_{lm}, \Gamma) + h)^2}. 
 \enm
\item  If \ $dist({\bf x}_{lm}, \Gamma)\le  \bar{W} h$,  then 
\eqm
 \left  | E^p_{lm} \right |  \le C |  \log h|. 
 \enm
\ei
\end{lemma}
 
 With the estimates in the above lemma, we prove the $L^2$ convergence of the pressure obtained from the three Poisson equations approach below.
\begin{theorem}
Assume that $p(x,y)\in C^2(\Omega^{\pm})$ and $h$ is small enough so that $dist(\Gamma,\partial \Omega)\sim O(1)$ and
the area enclosed by $A_{\Omega^-} \sim O(1)$, then the computed pressure from the three Poisson equation approach satisfies the following error estimate,
\eqm
  \|E^p\|_{L^2(\Omega_h)}\le C \sqrt{h} \, \left |\log h \right |.
\enm
Furthermore, 
\eqm
  \left | E^p_{lm} \right | \le C h, \qquad \mbox{\em if} \quad dist( {\bf x}_{ij}, \Gamma) \sim O(1). 
\enm
 \end{theorem}

 {\bf Proof: } Note that
 \eqmno
   \sum_{ij,dist({\bf x}_{ij}, \Gamma)>  \bar{W} h }  h^2 | E^p_{lm}  |^2 \le C \int_{0}^{2 \pi } d \theta \int_{h}^{\infty} \frac{r h^4}{(r + h)^4 }  dr \le C h .
\enmno
  From   the matrix-vector form of the  pressure error $\hat{A}_h \, {\bf E}^p  = {\bf T}^{p}$, we continue to estimate from~\eqn{eplm}
  \eqmno
  \|E^p\|^2_{L^2(\Omega_h)} &=& h^2 \sum_{ij} | E^p_{lm}  |^2 = h^2  \left (  \sum_{ij,dist({\bf x}_{ij}, \Gamma)>  \bar{W} h } +  \sum_{ij, dist({\bf x}_{ij}, \Gamma) \le \bar{W} h }  \right ) \\ \eqsp
  &\le &  C  h  +  h^2 N \bar{W} ( \log h)^2 \sim C h  ( \log h)^2,
\enmno
which leads to $\|E^p\|_{L^2(\Omega_h)}  \le  C\sqrt {h} \, |\!\log h|$.

\begin{remark}
 We think that the error estimates for the pressure that $\left |E^p_{lm} \right |\sim C |\log h| $ near the interface can be improved to $O(1)$; and $\left | E^p_{lm} \right |\sim C  h$ if $dis(\Gamma, {\bf x}_{lm})  \sim O(1)$ can be improved to 
 $dis(\Gamma, {\bf x}_{lm}) \le C \sqrt{h} $. 
\end{remark}

\subsection{Convergence of the velocity obtained from the three-Poisson equation approach}

Although we obtain the velocity by solving two Poisson equations, we can not use the  convergence theorem from \cite{li:mathcom} directly to obtain the  convergence since the finite difference approximation to $\grad p$ will generate $O(1/h)$ error. 
The convergence proof for the velocity is more challenging because of the involved gradient of the pressure. 
It is well-known that  computed velocity from the IB method has a first order convergence even though the pressure does not converge near the interface in the pointwise norm. There must be some cancellations around the interface.  

To describe the singularity due to the discretization of $p_x$ and $p_y$, we introduce the following discrete Dirac delta function that will be needed in the error estimates,
\eqm \label{delta-new}
  \delta^1_h(x) = \left \{ \begin{array}{ll}
 \dsp \frac{1}{h}    &  \mbox{if  $\dsp  |x| \le \frac{h}{2}$,}  \\ \eqsp
   0  &   \mbox{otherwise. }
  \end{array} \right.
\enm
This discrete delta function $\delta^1_h(x) \in L^2$ is probably the simplest, has the least support, and satisfies the zeroth consistency condition \cite{rpb-rjl:1d} as  that of the  hat and cosine discrete functions. We have applied the discrete delta function to  an example of Stokes equations in Section~\ref{sec-ex}. We can see that the numerical results are comparable to that obtained using the discrete cosine delta function. We use this discrete delta function in this paper mainly for theoretical purpose though.

From the jump conditions
\eqm
  \dsp  f_1 (s)  = [p] (s)  \cos \theta  - [u_n](s)  , \qquad  f_2(s)  = [p] (s)  \sin \theta  - [v_n](s),
\enm
 we know that the components of the velocity satisfy the following discrete Poisson equations
\eqmno
 \Delta_h u_{i,j}   &=&  \dsp  \frac{ P_{i+1,j}-P_{i-1,j} }{2h}  + \sum_{k=1}^{N_b} \left(   [p] (s_k) \cos \theta_k  -[u_n]_k  \frac{\null}{\null}\right  ) \delta_h(x_i-X_k)  \delta_h(y_j-Y_k) \Delta s_k , \\ \eqsp
 \Delta_h v_{i,j}   &=&  \dsp  \frac{P_{i,j+1}-P_{i,j-1}}{2h} + \sum_{k=1}^{N_b}  \left(   [p] (s_k) \sin \theta_k  -[v_n]_k \frac{\null}{\null} \right  )  \delta_h(x_i-X_k)  \delta_h(y_j-Y_k) \Delta s_k.
\enmno

Next, we consider the discrete relations that are analogous to the continuous situations.  These relations are needed in the velocity error estimate. The key point is that the singular sources can be decomposed as two parts: one part corresponds to the discontinuities in the gradient of the velocity; and other part corresponds to the discontinuity in the pressure. We start with the following lemma, which states that the singularity from the pressure is canceled out to some extent in the momentum equation discretization.

\begin{lemma} \label{lemma_px}
 Assume  we have a uniform grid $(x_i,y_j)$, $i,j=0,\cdots,N$. Let $p(x_{i+1},y_j)$ and $p(x_{i-1},y_j)$ be the pressure from different sides of the interface and assume that the interface intersects the grid line $y=y_j$ at $(X_{ij},y_j)$, then
 \eqml{pxdel2}
 && \dsp \lim_{h\goto 0}  \int \!\!  \int_{ dist({\bf x}, \Gamma)\le Wh }   p_x \, dx dy  = \int_{\Gamma} [p](s) \cos \theta(s) ds, \\ \eqsp  \eqsp
\null \!\!\!   && \dsp \sum_{ dist({\bf x}_{ij}, \Gamma)\le Wh }  h^2 \left (\frac{ p(x_{i+1},y_j)-p(x_{i-1},y_j)}{2h} \right )  =    \dsp \int_{\Gamma}  [p](s) \cos \theta(s) \,  ds + O(h).
  \enml

\ignore{
  \eqml{pxdel}
 \dsp   && \dsp \frac{ p(x_{i+1},y_j)-p(x_{i},y_j)}{h}  =  [p] \, \delta^1_h(x_i - X_i) ( {\bf n} \cdot {\bf e}_1) + O(1) \\ \eqsp
  &  & \qquad \null =  \dsp \int_c^d  [p]\,  \delta^1_h(x_i - X_i) ( {\bf n} \cdot {\bf e}_1) \delta^1_h(x_i - X_i) \delta(y-y_j) dy + O(1) ; 
  \enml
  \eqml{pxdel2}
 \dsp   && \dsp \sum_{ dist({\bf x}_{ij}, \Gamma)\le Wh }  \left (\frac{ p(x_{i+1},y_j)-p(x_{i},y_j)}{h} \right )  =    \dsp \int_{\Gamma}  [p] \cos \theta \,  ds,
  \enml
  where ${\bf e}_1$ is the unit vector in the $x$-direction. } 
\end{lemma}

\begin{figure}[hptb]
 \centerline{\includegraphics[width=0.45\textwidth]{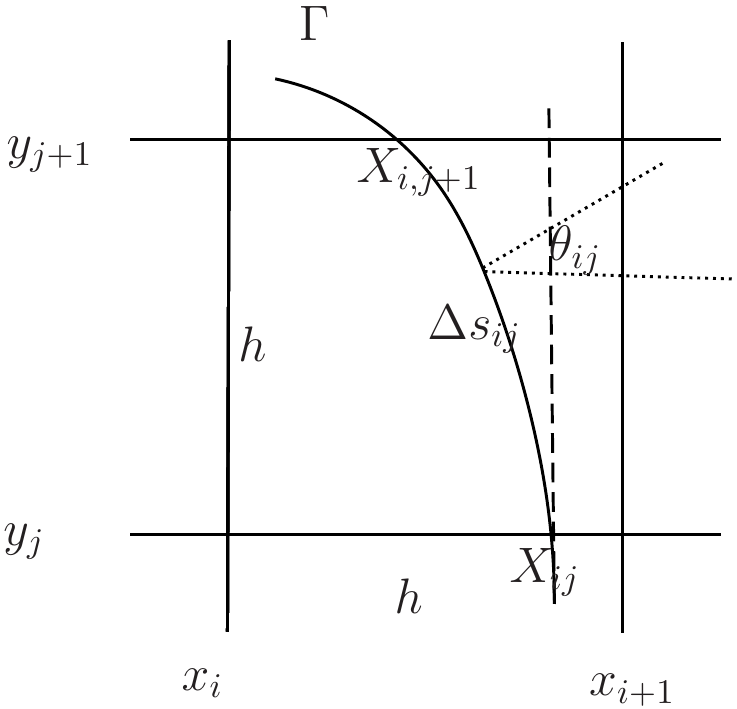}}
 \caption{A diagram of  an irregular grid point for the pressure where $(x_i,y_j)$ and $(x_{i+1},y_j)$ are from different sides of the interface. Note that  in this cell, we have  $\Delta y = h \approx \Delta s_{ij} \cos \theta_{ij}$.} \label{px-dy-ps}
\end{figure}

{\bf Proof:} \ Assume that $h$ is small enough so that the interface only cut the grid line $y_j$ only once between $x_{i-1}$ and $x_{i+1}$. In the continuous case, let $\Gamma_h$ be a closed domain containing the interface $\Gamma$, {\em cf.}, Figure~\ref{px-dy-ps}.  Then from the Green theorem, we derive
\eqmno
  \dsp   \int \!\!  \int_{ \Gamma_h }   p_x \, dx dy  =  \int \!\!  \int_{\Gamma_h}   \nabla p \cdot  {\bf e}_1 \, dx dy= \int_{ \partial \Gamma_h^+ }  p^+ {\bf e}_1 \cdot {\bf n} \,ds - \int_{ \partial\Gamma_h^- }  p^- {\bf e}_1 \cdot {\bf n} \,ds  +O(h) \approx  \int_{ \Gamma}   [p] \cos \theta \,  ds ,
\enmno
since $p$ is bounded, where ${\bf e}_1$ is the unit vector in the $x$-direction. Thus, as $\epsilon\goto h$, we have ${\bf e}_1 \cdot {\bf n} \goto \cos \theta$,
$\,\partial \Omega_h^{\pm} \goto \Gamma$, we have the first identity.

We define the first order extension of $p^-(x,y)$ from $ \Omega^-$ to  $\Omega^+$ as
\eqm
 p^-_e(x_{i+1},y_j)  = p^-(X_{ij} ,y_j) + p_x^-(X_{ij},y_j) (x_{i+1}- X_{ij} ).
\enm
The extension of $p^+(x,y)$ from $ \Omega^+$ to  $\Omega^-$ is defined similarly. Thus, we have
\eqmno
&& \frac{ p(x_{i+1},y_j)-p(x_{i-1},y_j)}{2h}  = \frac{ p^+ + p_x^+ \, (x_{i+1} - X_{ij} )   - p(x_{i-1},y_j) }{2h}   + O(h) \\ \eqsp
 &&  \null \qquad = [p]  \, \delta^1_h(x_i - X_{ij} ) + \frac{ p^-  + p_x^-  (x_{i+1} - X_{ij} )   + [p_x](x_{i+1} - X_{ij} ) -p(x_{i-1},y_j) }{2h}  + O(h)
 \\ \eqsp
 &&  \null \qquad = [p]  \, \delta^1_h(x_i - X_{ij} ) + \frac{p_e^-(x_{i+1},y_j) - p(x_{i-1},y_j) + [p_x](x_{i+1} - X_{ij} ) }{2h} +O(h) \\ \eqsp
 &&  \null \qquad = [p]  \, \delta^1_h(x_i - X_{ij} ) + p_x^-(x_{i-1} ,y_j) + O(h^2) + O(1) +O(h)  \\ \eqsp
 &&  \null \qquad = \dsp   [p] \, \mbox{sign}( ( {\bf n} \cdot {\bf e}_1)) \,  {\color{blue} \delta^1_h}(x_i - X_{ij} )  + O(1),  \\ \eqsp
  &&  \null \qquad = \dsp \int  [p]    \,  \delta^1_h(x_i - X_{ij})  \,  \delta(y-y_j) dy + O(1).
\enmno
Here we can see the use of the new discrete delta function $\delta^1_h(x)$. Since it is consistent, in the summary, it can be replaced by any other consistent discrete delta function  with an order $O(h)$ error.

Let $Q_{ij}$ be a grid function of any  two-dimensional function $Q(x,y)\in C(\Omega)$. From the discrete Green theorem we also know that 
\eqmno
 && \hspace{-0.8cm} \dsp \sum_{ij, p_x^{Irr}} h^2 \,  \frac{ p(x_{i+1},y_j)-p(x_{i-1},y_j)}{2 h}  Q_{ij}   \\ \eqsp
  && \null \qquad =   \sum_{ij, p_x^{Irr}} h^2  [p](X_{ij},y_j)  \cos \theta_{ij}\,  \delta^1_h(x_i - X_{ij})  \,  \delta^1_h(y_j-X_{ij } )  \Delta s_{ij}  Q_{ij}  + O(1)  \\ \eqsp
 && \null \qquad = \sum_{k=1}^{N_b} [p](X_k,Y_k)  \cos \theta_k\, Q(X_k,Y_k)\,  \Delta s_k + O(h) + O(1)   \\ \eqsp
  && \null \qquad = \int_{\Gamma}  [p](s)  \cos \theta (s) \, Q(X(s),Y(s)) \, d s + O(1).
  \enmno
Here we use the notation $\Delta s_{ij}$ because the intersections depend on both $i$ and $j$, and we have used the fact that $ \cos \theta_{ij}\, \Delta s_{ij} \approx  \Delta  y = h$,  see Figure~\ref{px-dy-ps} for an illustration.

Now we are ready to prove the first order convergence for the velocity. We use $u$ for the proof since the proof process is similar for the $v$
component.

\begin{theorem}
Let $u(x,y)$ be the solution to \eqn{Sa}-\eqn{Sc} and  $\,\mathbf{U}$ be the solution vector obtained from the
   immersed boundary method \eqn{IBa}-\eqn{IBc} using a discrete delta function. Then $\,\mathbf{U}$ is first order accurate with
   a logarithm factor in the infinity norm, that is,
\eqm
 |E_{ij}^u | \le C h | \log h|, \qquad i, j =1, 2, \cdots, N-1.
\enm
\end{theorem}

\textbf{Proof:}
Consider the error at a grid point $E_{lm}$,  if ${\bf x}_{lm}$ is close to the interface, that is,
$dist(\Gamma, {\bf x}_{lm})\le Wh$, we have
\eqml{eues}
  E_{lm}^u &=& \dsp \lp (A_h)^{-1} ( {\bf T}^{u}  -  {\bf T}^{p_x} ) \rp_{lm}   \\ \eqsp
 &=&  \lp (A_h)^{-1} ( {\bf T}_{reg}^{u}  -  {\bf T}_{reg}^{p_x} )  \rp_{lm} + \lp (A_h)^{-1} ( {\bf T}_{irr}^{u} -  {\bf T}_{irr}^{p_x} ) \rp_{lm}\\ \eqsp
  &=&  \dsp O(h^2) + \lp (A_h)^{-1}  ( {\bf T}_{irr}^{u}  -  {\bf T}_{irr}^{p_x} )  \rp_{lm} \\  \eqsp
  &=& \dsp \sum_{ dist({\bf x}_{ij}, \Gamma)\le Wh
  } \left\{  \lp h^2 T_{ij}^u  (A_h)^{-1} {\bf e}_{ij} \frac{1}{h^2} \rp_{lm}  -  \lp h^2 T_{ij}^{p_x}  (A_h)^{-1} {\bf e}_{ij} \frac{1}{h^2} \rp_{lm}  \right \} + O(h^2) \\ \eqsp
  &=& \dsp \sum_{ dist({\bf x}_{ij}, \Gamma)\le Wh }  h^2 \left\{ \frac{\null}{\null}   \lp \Delta_{h} u(x_i,y_j) - C_{ij}^{u}\rp   {\bf G}^{h}({\bf x}_{ij},{\bf x}_{lm})     \right.
  \\ \eqsp
  && \dsp \qquad  \qquad   \left.  \null - \lp D_h^x p (x_i,y_j)  -C_{ij}^{p_x}\rp
  {\bf G}^{h}({\bf x}_{ij},{\bf x}_{lm})   \right \} + O(h^2)  
\enml
where $C_{ij}^{u}$ and $C_{ij}^{p_x}$ are   from the source strength in the $x$- directions using
$f_1 (s)  = [p] (s)  \cos \theta  - [u_n](s)$,
\eqml{upcor}
  C_{ij}^{u} &=& \dsp \sum_{k=1}^{N_b} [u_n] _k \, \delta_h(x_i-X_k)  \delta_h(y_j-Y_k) \Delta s_k = \int_{\Gamma} [u_n](s) \,
  \delta_2({\bf X}_{ij} - {\bf X} (s)) ds + O(h) \\ \eqsp
  C_{ij}^{p} &=& \dsp \sum_{}^{} [p] ({\bf X}_{ij} )  \cos \theta_{ij} \, \delta_h(x_i-X_{ij})  \delta_h(y_j-Y_{ij}) \Delta s_{ij}  = \int_{\Gamma} [p]  \cos \theta \,
  \delta_2({\bf X}_{ij} - {\bf X} ) ds + O(h).
\enml
Thus, the error is consists of two parts: one if $p_x$ that  is offset part of singular  source term, which we have just discussed; the second part is due to the source distribution that will offset other parts of singular source term up to $O(h)$, which we can apply the estimates of the IB method for the scalar Poisson equations with an immersed closed interface. 
From Theorem~3.8 in \cite{li:mathcom}, we know that the first term in the last line of \eqn{eues} is bounded by 
$C h |\log h|$.
For the second term, from Lemma~\ref{lemma_px}, we know that
\eqmno
  &&  \hspace{-0.8cm} \dsp \sum_{ dist({\bf x}_{ij}, \Gamma)\le Wh } \!\!\! \lp h^2   D_h^x p (x_i,y_j) {\bf G}^{h}({\bf x}_{ij},{\bf x}_{lm})  -C_{ij}^{p_x}    {\bf G}^{h}({\bf x}_{ij},{\bf x}_{lm})   \rp  \\ \eqsp
  && \null \qquad =  \sum_{ij, p_x^{Irr}}    h^2   D_h^x p (x_i,y_j) {\bf G}^{h}({\bf x}_{ij},{\bf x}_{lm})   -  \!\!\!   \!\!\! \sum_{ dist({\bf x}_{ij}, \Gamma)\le Wh }  \!\!\!   h^2    C_{ij}^{p_x}    {\bf G}^{h}({\bf x}_{ij},{\bf x}_{lm}) \\ \eqsp
  && \null \qquad =  \int_{\Gamma}  [p](s)  \cos \theta (s) \, {\bf G}_I^{h}({\bf X}(s),{\bf x}_{lm}) \, d s + O(h) - \!\!\!  \!\!\! \sum_{ dist({\bf x}_{ij}, \Gamma)\le Wh }  \!\!\!  \!\!\!   h^2    C_{ij}^{p_x}    {\bf G}^{h}({\bf x}_{ij},{\bf x}_{lm})  \\ \eqsp
  && \null \qquad = O(h).
\enmno

If $dist(\Gamma, {\bf x}_{lm})> Wh$, the proof above is still valid except that we are not going to have the singular integration. Thus, we do not need  to have the $|\log h|$ factor. This means that for the IB method, the larger errors often occur near or on the interface. This completes the proof.

\ignore{
\begin{remark}
  For convenience of the proof, we assume that the local truncation error of the treatment of the pressure boundary is $O(1)$ so that the error can be $O(h)$ near the boundary. However, for the $\sqrt{h} | \log h |$ convergence in the $L^2$ norm, the condition can be relaxed to $O(1/h)$ for the local truncation error at the numerical boundary points for the pressure.
\end{remark} }

\section{Convergence analysis of  the IB method for Stokes equations using a MAC grid} \label{sec:MAC}

It is more challenging to prove the convergence of the pressure and velocity obtained from the IB method using a MAC grid because the pressure and the velocity are couple together.  Our strategy is,  similar to the continuous case, to write the linear system of equations as two parts. The first part is for the pressure equations,  which  does not contain the velocity from the discrete divergence operation at interior points.  At boundary points, we can not get rid of the velocity but we can get the estimates that is good enough for the proof of the pressure. The second part of the linear system of equations is for the velocity equations, which has similar representations as that of the three Poisson equation approach. The differences is the way of discretizing the pressure and at different grid points for the primary variables.

\ignore{
The discussion for 2D problems is much more challenging since the interface is often a curve instead of a point.
In \cite{mori-proof}, the author has proved the first order convergence of the IB method for the Stokes equations with a periodic boundary condition in 2D based on existing estimates
between the discrete Green function and the continuous one in~\cite{Hasimoto-stokes56}.
However, there are almost no theoretical proofs on the IB method for elliptic interface problems or other PDEs with general boundary conditions. We will prove that the result obtained from the IB for the elliptic interface problem with a Dirichlet boundary condition is indeed first order accurate in this section.
}

\ignore{
Consider the following two dimensional (2D) Stokes equations,
\eqm
 && \nabla p = \mu \Delta  {\bf u}  + {\bf G} + \int_{\Gamma} {\bf f}(s)\,
\delta_2({\bf x}-{\bf X}(s)) ds, \quad {\bf x} \in \Omega, \label{Sa}\\ \eqsp
&& \nabla \cdot {\bf u} = {\bf 0}, \quad {\bf x} \in \Omega, \label{Sb} \\ \eqsp
&& \left. {\bf u} (x,y)\right |_{\partial \Omega}={\bf u} _{0}(x,y)
\enm
where we assume that ${\bf x}=(x,y)$, ${\bf u} =(u,v)$,  $\mu$ is a constant, $s$ is a parameter, {\em e.g.}, the arc-length of the interface $\Gamma$,
${\bf G} \in C(\Omega)$, $\Gamma(s)\in C^1$, $ {\bf f}(s)\in C^1$, see Fig.~\ref{domainD} for an illustration.}


\subsection{The IB method for solving Stokes equations on a MAC grid}

\begin{figure}[hptb]
 \centerline{  \includegraphics[width=0.75\textwidth]{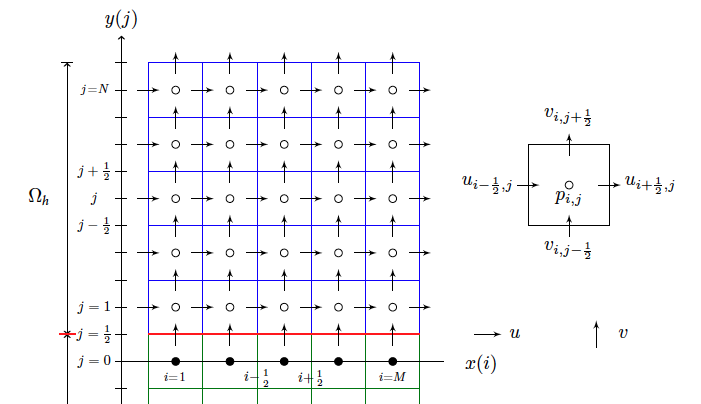} }
 \caption{A diagram of a MAC grid adapted from \cite{cheng-Kian-Lai-18}.} \label{MAC-grid}
\end{figure}

We adopt the notations from \cite{cheng-Kian-Lai-18}, see also Figure~\ref{MAC-grid} for an illustration,  by setting a MAC grid,
\eqm
 x_i = a + \left ( i+\half \right ) h, \quad i=0,1,\cdots, N;  \qquad y_j= c   +  \left (j+\half \right ) h, \quad j=0,1,\cdots, N.
\enm   
The discrete form of the momentum equation for $u$ is
\eqm \label{IBa}
 \mu \frac{U_{i-1/2,j}+U_{i+3/2,j}+U_{i+1/2,j-1}+U_{i+1/2,j+1}-4U_{i+1/2,j}}{h^2} - \frac{P_{i+1,j}-
P_{i,j}}{h} = \dsp F^{1}_{i+1/2,j},
\enm
for $i=1,2,\cdots, N-1$, $j=1,2,\cdots, N$. There are $N(N-1)$ equations.  In the representation above,  for example, the force component
$F^{1}_{i+1/2,j}$ is defined as
\eqml{Fij}
F^{1}_{i+1/2,j} &=&  \dsp G^{1}_{i+1/2,j}  + \sum_{k=1}^{N_b} f^1_k \, \delta_h(x_{i+1/2}-X_k) \, \delta_h(y_j-Y_k) \, \Delta s_k \\ \eqsp
  &=&  \dsp G^{1}_{i+1/2,j}  + \sum_{k=1}^{N_b} \left ( \hat{f}^1_k \, n^1_k - \hat{f}^2_k n^2_k \right ) \delta_h(x_{i+1/2}-X_k) \, \delta_h(y_j-Y_k) \, \Delta s_k,
\enml
with ${\bf n}_k = ( n^1_k, n^2_k)$ being the unit normal direction at $(X_k,Y_k)$, where $\delta_h(x)$ is a discrete delta function such as a hat,  or cosine, or radial discrete delta function. 

Similarly the discrete form of the momentum equation for $v$ is
\eqm \label{IBb}
\mu \frac{V_{i-1,j+1/2}+V_{i+1,j+1/2}+V_{i,j-1/2}+V_{i,j+3/2}-4V_{i,j+1/2}}{h^2} - \frac{P_{i,j+1}-
P_{i,j}}{h} = \dsp F^{2}_{i,j+1/2},
\enm
for $i=1,2,\cdots, N$, $j=1,2,\cdots, N-1$. There are $(N-1)N$ equations. We also discretize the divergence free condition according to
\eqm
\frac{U_{i+1/2,j}-U_{i-1/2,j} }{h} + \frac{V_{i,j+1/2}-V_{i,j-1/2}}{h} = \dsp 0, \label{IBc}
\enm
for $i=1, \cdots, N, \quad j=1,\cdots,  N$.
The total number of equations is $2N(N-1)+N^2$.  The total unknowns for $U_{ij}$ is $(N-1)N$, so is for  $V_{ij}$, and the total unknowns for $P_{ij}$ is $N^2$. Thus, the number of equations is the same as the number of unknowns. If we can solve the above linear system of equations, then we get an approximate solution for the velocity and the pressure that can  differ by a constant. There is a large collection of references in the literature on how to solve the linear system, which is a saddle point problem.   Our focus in this paper is to obtain error estimates  for the discrete solutions of  the discrete system of equations.

\subsubsection*{A sub-linear system of equations}

It is well-known that the discrete system of equations \eqn{IBa}, \eqn{IBb}, and \eqn{IBc} has a unique set of solutions for the velocity, and the pressure if we fix the pressure if we fix the pressure  at one point, say $P_{11}=0$. We denote the solution by $U^*_{i+1/2,j}$, $V^*_{i,j+1/2}$, and $P_{ij}^*$ and try to find error estimates for the  approximations to the  Stokes equations. We use the $*$ to represent the solution to the discrete system to distinguish from the discretization in the hope to make the proof clearer.  

We divide by $h$ from the equation \eqn{IBa} and  \eqn{IBb}  at an interior grid point $(x_i,y_j)$, and divide by $-h$ from the equation \eqn{IBa} and  \eqn{IBb}  at $(x_{i-1},y_{j-1})$ for $i,\, j =2,\cdots, N-2$,  and then add them together. Applying  equation \eqn{IBc} with  some tedious manipulations, we  get
\eqml{plap}
 &  \dsp  \frac{P^*_{i-1,j}+P^*_{i+1,j}+P^*_{i,j-1}+P^*_{i,j-1 }-4P^*_{i,j }}{h^2}  =\frac{F^{1}_{i+3/2,j} - F^{1}_{i+1/2,j}}{h} +  \frac{F^{2}_{i,j+3/2} - F^{2}_{i,j+1/2}}{h}, \\ \eqsp
 & i=2, \cdots, N-2, \qquad j=2,\cdots, N-2.
\enml
The process consists in elementary row transformations applied to the system of equations \eqn{IBa}-\eqn{IBb}  from the point of view  of linear algebra.  Thus, the solution remains unchanged. 
If we apply the divergence operator to the moment and apply the divergence theorem,  we would get a Poisson equation for the pressure without  the velocity. We would get the above discrete system if we apply the standard five-point discrete Laplacian. The complication would be the boundary condition. 
The discrete momentum equation at $i,j=1$,  $i,j =N-1$ can be regarded as a Neumann boundary condition for the pressure, even though they are still coupled with other equations,  for example,
%
\eqm \label{pbc_a}
  \frac{P^*_{2,j}-P^*_{1,j}}{h} = \mu \frac{U^*_{1/2 ,j}+U^*_{5/2,j}+U^*_{3/2,j-1}+U^*_{3/2,j+1}-4U^*_{3/2,j}}{h^2} -  F^{1}_{3/2,j},
\enm
$j=1,\cdots,N-1$.
In other words, the approximate pressure is the solution of the discrete Laplacian with an approximate Neumann boundary condition that is coupled with the velocity near the `boundary' since no pressure boundary condition is needed or used. 
If we can estimate the right hand side of the above, as an approximation of the artificial boundary condition, then we can get an estimate for the pressure using previous analysis and results.

\ignore{
along with the boundary condition discretization that can be written as
\eqm
  B_h({\bf P}_h, {\bf U}_h)=0,
\enm
where ${\bf P}_h, {\bf u}_h$ are the grid values of $p$ and ${\bf u}$, respectively.
The above system of equations is a system of equations of the discrete Laplacian for the pressure. }

In fact, the  normal derivative of the pressure  along the boundary  is not a free variable and it should satisfy
\eqm
  \frac{\partial p}{\partial n} = \mu \Delta {\bf u}  \cdot {\bf n} + {\bf G} \cdot {\bf n}, \qquad {\bf x} \in \partial \Omega, \label{real-p-bc}
\enm
from the momentum equation.  The above relation is true for any closed region within the solution domain. Note that for a rectangular domain, we know that ${\bf n} = (\pm 1, 0)$ or $(0,\pm 1)$,  and the discrete equations
in \eqn{IBa}  and \eqn{IBb}  are the discrete form of \eqn{real-p-bc} when $i,j=1$ or  $i,j=N-1$.
The finite difference  equation in \eqn{pbc_a} is a discretization of the above boundary condition at $(x_{3/2},y_j)$ corresponding to $i=1$. For a no-slip boundary condition, we have the  estimates for the discretization in the $x$-direction in the following lemma. 
 Similar estimates can be easily generalized to the $y$-direction as well.

\begin{lemma}  \label{p-trunc} Let the local truncation error of the discrete moment equation for $u$ be
\eqm
 T_{BC}^{p,x} = \frac{p(x_2, y_j )  - p(x_1,y_j ) }{h} - \mu \Delta_h u(x_{3/2},y_j)   -  F^{1}(x_{3/2},y_j) .
\enm
Then we have $|T_{BC}^{p,x} | \le C h^2$,  $|\Delta_h u(x_{3/2},y_j) |\le C$, and \ $\dsp \frac{p(x_2, y_j )  - p(x_1,y_j ) }{h} \sim O(h)$. 
\end{lemma}

{\bf Proof:} The first estimate is obtained using the Taylor expansion at $(x_{3/2},y_j)$ since the finite difference discretizations are centered second order accurate  ones for the $p_x(x_{3/2},y_j)$ and $\Delta u(x_{3/2},y_j)$.

For the second part, $\Delta_h u(x_{3/2},y_j)$  consists of  two parts, that is, the second order derivative in the $x$ and $y$ directions. For the approximation of $u_{yy}(x_{3/2},y_j)$, we know that
\eqmno
   && \frac{u(x_{3/2},y_{j-1})  - 2u(x_{3/2},y_{j} ) + u(x_{3/2}, y_{j+1} )   }{h^2} \\ \eqsp
   && \qquad \hspace{2cm}\null  =u_{yy}(x_{3/2},y_j) + O(h^2) = u_{yy}(x_{1/2},y_j) +O(h) =0 +  O(h)
\enmno
since $u(x_{1/2},y)=0$ and its tangential derivatives are also zero.

Now we check the second term, the finite difference approximation for $u_{xx}(x_{3/2},y_j)$. Note that since $u_x + v_y=0$ and $v=0$ along the boundary $x=x_{1/2}$, we have $v_y(x_{1/2}, y)=0$, so is
$u_x(x_{1/2}, y)=0$. Thus, we have
\eqml{uxx_exp}
  u(x,y) &= & \dsp u(x_{1/2},y) + (x-x_{1/2}) u_x (x_{1/2},y)  +\frac{(x-x_{1/2})^2}{2}   u_{xx} (x_{1/2},y) + O(h^3)\\ \eqsp
    &= & \dsp \frac{(x-x_{1/2})^2}{2}  u_{xx} (x_{1/2},y)  + O(h^3).
\enml
Therefore, the second  term in $\Delta_h u(x_{3/2},y_j)$   satisfies
\eqml{uxx_exp2}
   && \dsp \frac{u(x_{1/2},y_{j})  - 2u(x_{3/2},y_{j} ) + u(x_{5/2},y_j)   }{h^2}  \\ \eqsp
   &&  \hspace{1cm} \dsp  \null = \frac{-2 u_{xx} (x_{1/2},y_j) \frac{ h^2}{2}  +  u_{xx} (x_{1/2},y_j) \frac{ 4 h^2 }{2}  + O(h^3) }{h^2}
   = u_{xx} (x_{1/2},y_j) + O(h) \\ \eqsp
   &&  \hspace{1cm} \dsp  \null = \, O(1).
\enml

This lemma also explains why the pressure is only first order accurate when we use $\frac{\partial p}{\partial n}=0$ as an approximate pressure boundary condition, whose discrete form in a rectangular domain can be written as
\eqm \label{pbc-0}
  \frac{P_{2,j}-P_{1,j}}{h}=0, \quad  \frac{P_{N,j}-P_{N-1,j}}{h}=0, \quad  \frac{P_{i,2}-P_{i,1}}{h}=0, \quad  \frac{P_{i,N}-P_{i,N-1}}{h}=0.
\enm
More accurate but  sophisticated treatment of numerical pressure boundary conditions can be found in \cite{Johnston-Liu02} for rectangular domains, and \cite{xwan-li-pbc} for more
general geometries. In our error analysis, for the $L^2$ convergence of the pressure, we just need the local truncation error being bounded by $O(1/h)$ at the `boundary grid points'.  For a first order pointwise convergence, in general we need the local truncation error that can be  bounded by $O(1)$, and for second order point-wise convergence, we need the local truncation error  that can be  bounded by $O(h)$.

Since the velocity is solved from the momentum equation, the proof of the velocity is almost the same as that of the three Poisson equation approach. The grid points location differences contribute at most $O(h)$ difference that have no effect on the first order convergence of the velocity. 

\section{A numerical example to verify the discrete delta function and convergence behaviors }  \label{sec-ex}

The main purpose of the numerical results is to verify our analytic analysis. As a matter of fact, numerical experiments  have helped and guided our theoretical analysis.

Now we consider an example with a non-smooth velocity.
The exact solution is the following:
\eqml{example}
u(x,y) &=&\left \{ \begin{array}{ll}
 \dsp  \frac{y}{2} \left (x^4  - y^4 + 2 y^2 -1\right )   & \mbox{if  $r>1$,   } \\ \eqsp
\dsp \frac{y}{2} \left (x^2  + y^2 -1\right )     &    \mbox{ if $r \le 1$,}
\end{array} \right.  \\ \eqsp
v(x,y) &=&\left \{ \begin{array}{ll}
\dsp - x^3 \left (x^2  + y^2 -1\right )   & \mbox{if  $r>1$,   } \\ \eqsp
\dsp - \frac{x}{2} \left (x^2  + y^2 -1\right )     &    \mbox{if $r \le 1 $,}
  \end{array} \right.  \\ \eqsp
p(x,y) &=& \left \{ \begin{array}{ll}
\dsp   -\frac{1}{2} x y   & \mbox{if  $r>1$,   } \\ \eqsp
\dsp  \frac{1}{2} x y    &    \mbox{ if $r \le 1 $,}
  \end{array} \right.
\enml
where $r = \sqrt{x^2 + y^2}$. The interface is the circle $r=1$
and the solution domain is $[-2,\; 2]\times [-2,\; 2]$.
Note that  this is a non-trivial example since  all the jump conditions are functions of $x$ and $y$  with a curved interface.

\begin{table}[hptb]
\caption{A grid refinement analysis of the IB method for \eqn{example} using the discrete cosine delta function. }
 \label{tab-ib-cos}

\vthin

\begin{center}
\begin{tabular}{|c| c c| c c|cc| }
  \hline
       $N$   & $\|E_{\bf u} \|_{\infty} $  & $r$  & $\| E_{ p} \|_{L^2} $  & $r$   & $\| E_{ p,\sqrt{h} }\|_{\infty} $  & $r$    \\    \hline \hline
  32 &     4.4971e-02    &        &   8.4473e-02      &      &   1.2208e-01       &      \\
 64 &   1.0813e-02 &  2.0562  & 2.8834e-02  & 1.5507   & 2.0772e-02  &  2.5551\\
 128&   2.7984e-03   & 1.9501 &  1.9107e-02  &  0.5936  &  1.2731e-02 &  7.0633 \\
 256 &  9.0206e-04 &  1.6333  & 1.3151e-02   & 0.5389  &  5.9490e-03   & 1.0976 \\
  512&  2.0304e-04   & 2.1515  &  9.0486e-03 &  0.5.395   &  1.8596e-03   & 1.6777   \\   \hline
Average & & 1.9478  & & 0.8057 & & 1.5092 \\ \hline \hline
\end{tabular}

\end{center}
\end{table}

In Table~\ref{tab-ib-cos} and \ref{tab-ib-new}, we show the grid refinement results using the discrete cosine and the new discrete delta function defined in \eqn{delta-new}, respectively. The second column is the error of the total velocity defined as
\eqm
  \|E_{\bf u} \|_{\infty} = \sqrt{ \frac{\null}{\null}  \|E_{ u} \|_{\infty} ^{2 \frac{\null}{\null} }   + \|E_{ v} \|_{\infty} ^{2 \frac{\null}{\null} }  }.
\enm
We can see that the discrete cosine delta function performs better than the simplest one \eqn{delta-new}. All the convergence orders have the expected or better results, that is, first order in the infinity norm for the velocity, half order $O(h^{1/2}$) for the pressure, and first order in the infinity norm for the pressure measured $\sqrt{h}$ away from the interface.

\begin{table}[hptb]
\caption{A grid refinement analysis of the IB method for \eqn{example} using the new discrete  delta function $\delta_h^1(x)$. }
 \label{tab-ib-new}

\vthin

\begin{center}
\begin{tabular}{|c| c c| c c|cc| }
  \hline
       $N$   & $\|E_{\bf u} \|_{\infty} $  & $r$  & $\| E_{ p} \|_{L^2} $  & $r$   & $\| E_{ p,\sqrt{h} }\|_{\infty} $  & $r$    \\    \hline \hline

  32 &  1.3149e-01    &           & 9.2928e-02       &        & 1.2852e-01      &       \\
  64 &   4.1643e-02   & 1.6588   & 3.0121e-02  & 1.6254  &  2.3241e-02  & 2.4672 \\
 128 & 1.5063e-02   & 1.4671   & 2.3433e-02  & 0.3622 &  1.4059e-02  & 0.7251\\
  256 & 6.6595e-03   & 1.1775   & 1.5072e-02 &  0.6367 &  1.1841e-02  & 0.2477  \\
  512 &  2.8297e-03   & 1.2348   & 1.1777e-02 &  0.3558  &  2.9406e-03  & 2.0096\\    \hline
Average & & 1.1076   & & 0.5960 & & 1.0899 \\ \hline \hline
\end{tabular}

\end{center}
\end{table}

\section{Conclusions, discussions, and acknowledgments}

 In this paper, we have shown the well-known fact that the immersed boundary method is first order accurate for the velocity in the pointwise (infinity) norm for Stokes equations with an immersed closed interface and a Dirichlet boundary condition. We also have proved that the pressure has half-order convergence $O(\sqrt{h}|\log h |$) in the $L^2$ norm, and first order convergence  ($O(h$)  away  from the interface  ($O(1)$). The key in the proof is first to separate the pressure from the velocity, and then use the discrete Green functions to get the estimate. We think it is possible to improve the estimate for the pressure that is $O(h$ accurate away  from the interface  ($O(\sqrt{h})$) instead of $O(1)$, which has been observed numerically. 
For the velocity, the singularity from the $\nabla p$ across the interface is offset by part of the singular source term, which can be considered as a source distribution of a special discrete delta function in $L^2$.

 While the proof is for Dirichlet boundary conditions for the velocity, we believe that the conclusions  are valid for other 
 boundary conditions  as long as the problems are wellposed since dominant errors are from the singular source and its discretization. From the proof, we also can see that for a first order method, we can use three Poisson equations for the pressure and the velocity without affecting the accuracy although the divergence free condition may not be satisfied exactly but to $O(h^2)$.

We would like to thank Drs. Ming-Chih Lai, Qinghai Zhang, and Hongkai Zhao for valuable discussions and suggestions.
Z. Li is partially supported by a Simons grant (No. 633724).

\bibliographystyle{plain}
\bibliography{../TEX/BIB/bib,../TEX/BIB/other,../TEX/BIB/zhilin,Mac_Stokes}

\appendix


\ignore{

Constructing Neumann Green functions.

/Users/zhilin/Dropbox/Research/IBM_Stokes/Neumann_Green/Green_Sum

}

\end{document}